\RequirePackage{fix-cm}
\documentclass[twocolumn]{svjour3}          
\smartqed  
\usepackage{graphicx}
\usepackage[utf8]{inputenc}
\usepackage[T1]{fontenc}
\usepackage{amssymb,amsmath,setspace}
\usepackage{bm}
\usepackage{subfig}
\usepackage{tikz}
\usepackage{algpseudocode}
\usepackage{algorithm}
\usepackage[%
  colorlinks%
  ,bookmarks={false}%
  ,pdftitle={}%
  ,pdfsubject={}%
  ,pdfkeywords={}%
  ,pdfauthor={emma.m.naden@gmail.com, maerz@maths.ox.ac.uk}%
  ]{hyperref}

\newcommand{\Real}{\mathbb{R}}
\newcommand{\surf}{\mathcal{S}}
\newcommand{\len}{\mathcal{L}}
\newcommand{\area}{\mathcal{A}}
\newcommand{\band}{B(\surf)}

\DeclareMathOperator{\Div}{div}
\DeclareMathOperator{\cp}{cp}

\DeclareMathOperator{\lingen}{span}

\newcommand{\pd}{\partial} 
\newcommand{\bd}{\partial} 
 
\newcommand{\T}{\text{T}}

\newcommand{\lap}{\Delta}
\newcommand{\grad}{\nabla}
\newcommand{\slap}{\lap_{\surf}}
\newcommand{\sgrad}{\grad_{\surf}}
\newcommand{\sDiv}{\Div_{\surf}}
\newcommand{\sInt}{\int\limits_{\surf}}
\newcommand{\Tang}[1]{{\mathcal{T}_{#1} \surf}}
\newcommand{\rel}{\text{rel}}

\newcommand{\dx}{h}
\newcommand{\dt}{\tau}
\newcommand{\ds}{\,\,.\!\!*}

\newtheorem{principle}[theorem]{Principle}

\graphicspath{ {./Figures/} }

\journalname{In preparation}
\begin{document}

\title{Anisotropic Diffusion Filtering of Images on Curved Surfaces}

\titlerunning{Anisotropic Diffusion on Curved Surfaces}

\author{Emma Naden \and Thomas M\"arz \and Colin B. Macdonald}

\authorrunning{E. Naden, T. M\"arz, C.B. Macdonald} 

\institute{
	Emma Naden \at
	Leidos Health\\
	705 E. Main St.\\
	Westfield, IN 46074, USA\\
	\email{emma.naden@oxon.org}           
	\and
	Thomas M\"arz \at
	Mathematical Institute\\
	University of Oxford\\
	Oxford OX2 6GG, UK \\
	\email{maerz@maths.ox.ac.uk} 
	\and
	Colin B. Macdonald \at
	Mathematical Institute\\
	University of Oxford\\
	Oxford OX2 6GG, UK \\
	\email{macdonald@maths.ox.ac.uk} 
}

\date{Received: \today / Accepted: date}

\maketitle

\begin{abstract}
We demonstrate a method for filtering images defined on curved surfaces embedded in 3D.
Applications are noise removal and the creation of artistic effects. 
Our approach relies on in-surface diffusion:
we formulate Weickert's edge/coherence enhancing diffusion models in a surface-intrinsic way.
These diffusion processes are anisotropic and the equations depend non-linearly on the data. 
The surface-intrinsic equations are dealt with the closest point method, a technique for solving
partial differential equations (PDEs) on general surfaces.
The resulting algorithm has a very simple structure: we merely
alternate a time step of a 3D analog of the in-surface PDE
in a narrow 3D band containing the surface with a reconstruction of the surface function.
Surfaces are represented by a closest point function.
This representation is flexible and the method can treat very general surfaces.
Experimental results include image filtering on smooth surfaces,
open surfaces, and general triangulated surfaces.

\end{abstract}

\section{Introduction}\label{sec:intro}
Image processing is an active field of mathematical research,
and models based on partial differential equations (PDEs) have been
used successfully for various image processing tasks.
Examples include inpainting
\cite{bscb00,sc02,m11},
segmentation \cite{cv01,hkiu07,zzsw11},
smoothing and denoising
\cite{pm90,w98,clmc92,alm92}.

Denoising on general curved surfaces \cite{lai2011framework}, the
corresponding scale space analysis \cite{Kimmel97intrinsicscale}, and
related image processing problems have seen some interest
\cite{spira,bogdanova}.
Of particular relevance for our voxel-based technique is that image
processing on curved surfaces can occur even when data is acquired in
three-dimensional volumes.
For example, Lin and Seales \cite{Lin:2005} propose CT scanning of
scrolls in order to non-destructively read text written on rolled up
documents.
Other applications of denoising on surfaces include
digital image-based elasto-tomography (DIET), a technique for
non-invasive breast cancer screening \cite{brown2010vision},
texture processing \cite{clarenz2004processing,bajaj2003anisotropic},
and surface fairing where the surface data is itself a height-field
perturbation relative to a reference surface
\cite{desbrun2000anisotropic}.

In the cited articles a common theme is the use of PDE models.
The topic of solving PDEs on curved surfaces is an important area of research with many additional applications in physics and biology
(for example to model pattern formation via reaction-diffusion on
animal coats \cite{m81,sw03} or the diffusion of chemicals on the
surface of a cell \cite{sacao05,ngcrss07,f02}).

The numerical treatment of PDEs on surfaces requires a suitable representation of the geometry.
The choice of the representation is fundamental to the complexity of the algorithm.
Parametrizations are commonly used, however they can introduce artificial distortions and singularities into the models
even if the surface geometry is as simple as that of a sphere.
Moreover, joining up multiple patches is typically necessary if the geometry is more complex.
On triangular meshes approximating a surface, discretizations can be obtained using finite elements, e.g., \cite{Dziuk:2007}.
Embedding techniques based on implicit representations such as level sets (e.g., \cite{bertalmio}) or the closest point representation
are further alternatives and tend to be quite flexible.
However with level sets it is not obvious how to treat open surfaces
and one has to introduce artificial boundary conditions at the boundary of
the 3D neighborhood in the embedding space.

In this work we use the closest point method \cite{Ruuth/Merriman:jcp08:CPM} for solving the diffusion models on surfaces.
This technique keeps the resulting evolution as simple as
possible by alternating between two straightforward steps:
\begin{enumerate}
\item A time step of the analogous three-dimensional PDE model
  using standard finite difference methods.
\item An interpolation step, which reconstructs the surface function and 
  makes the 3D calculation consistent with the surface problem.
\end{enumerate}
The surface geometry is encoded in a closest point function \cite{Ruuth/Merriman:jcp08:CPM,Maerz/Macdonald:cpfunctions}
which maps every off-surface point to that surface point which is closest in Euclidean distance.
The closest point function is used in the reconstruction step only and therefore this technique does not require modification
of the model via a parameterization nor does it require the surface to
be closed or orientable.
The closest point method has already been used successfully in
image segmentation \cite{luke:segment}, visual effects \cite{Auer/Macdonald/Treib/Schneider/Westermann:fluidEffects}
and to perform Perona--Malik edge-stopping diffusion \cite{bvgmm13}.

The use of diffusion processes to smooth noisy images dates back to at least the 1980s, e.g.,
Koenderink \cite{k84} notes that one could perform smoothing by solving the heat equation over the image domain. 
Unfortunately, the linearity of this equation implies a uniform blur of the entire image without any consideration of the structure of the image.
Consequently, important high-frequency components due to edges are blurred as much as 
undesired high-frequency components due to noise and the structure of the image is lost. 
Perona and Malik \cite{pm90} put forth an elegant and simple modification of the heat equation to address this issue.
In order to avoid blurring edges, they vary the rate of diffusivity according to the local gradient of the image,
stopping the diffusion process at edges (i.e., where the magnitude of the gradient is large).
Regularizations of the Perona--Malik model and other variations of the diffusivity 
(and also relations to geometric evolutions such as mean curvature motion) have been studied for example in \cite{clmc92,alm92}.
Although there is an improvement over linear diffusion, Perona and Malik's method retains a fundamental flaw:
because the diffusion coefficient drops to zero when the gradient is high, their technique cannot remove noise along edges.
This motivates the use of an anisotropic diffusion filter as suggested by Weickert \cite{w98}
in which one alters not only the rate of diffusion near an edge but also its direction 
so that diffusion is performed along rather than across the edge.
In order to accomplish this, the scalar diffusion coefficient found in linear and nonlinear Perona--Malik diffusion 
must be replaced with a diffusion tensor.
Apart from edge-enhancing denoising, Weickert's diffusion model \cite{w98} has a second application called coherence-enhancing diffusion
which can be used to create artistic effects.
The contribution of this paper is to extend Weickert's model to curved surfaces and 
implement anisotropic diffusion on surfaces using the closest point method.
Anisotropic diffusion to smooth surface meshes as well as functions defined on surface meshes has been considered in \cite{bajaj2003anisotropic}.
Their approach differs from ours in two main ways.
Firstly, the design of their diffusion tensor is based on the directions of principal curvature
while ours is based on eigenvectors of a surface-intrinsic structure tensor.
Secondly, they use surface finite elements to discretize the diffusion equation
while we use the closest point method.

The rest of this paper unfolds as follows.
Section~\ref{sec:image_processing} reviews the concepts of anisotropic diffusion filtering closely following Weickert \cite{w98}.
In Section~\ref{sect:aniso_surf} we introduce a mathematical model for surface-intrinsic aniso\-tropic diffusion.
In particular, we define the structure tensor for images on surfaces which allows us to estimate visual edges.
Section~\ref{sec:cpm} introduces the closest point method for solving diffusion PDEs posed on surfaces
and describes our numerical schemes. 
The presentation of results and a discussion follows in Section~\ref{sec:results}.

\section{Anisotropic Diffusion in Image Processing} \label{sec:image_processing}
A \emph{digital image} is a quantitative representation of the information (color and intensity values) encoded in an image.
Typically, the domain of a digital image is a finite and discrete set of points called \emph{pixels}.
We view images as a functions defined on a rectangular domain in $\Real^2$.
Digital images are then discretizations of such functions obtained from sampling on a Cartesian grid. 
We will later extend this to images defined on two-dimensional surfaces.   
 
When we view images as a functions, many image processing tasks can be described by time-dependent PDE processes such as diffusion.
The resulting image is then the solution to the PDE. 
Color images are represented by vector-valued functions, for example, in
RGB color space we have a triple $u(x) = [R,G,B]$ at each point $x$ in the domain, where $R,G,B \in [0,1]$ represent the intensities of red, green, and blue respectively.
When working with color images, we apply the diffusion filter to each of the color components individually.  Unless otherwise noted,
throughout the remainder of this paper we assume $u \in [0,1]$, where zero represents black (no color intensity) and one represents white (full color intensity).

A so-called \emph{noisy image} is an image that differs from the real object depicted (the \emph{ground truth}),
by a discoloration of pixels in the digital image. Noise can occur from a variety of sources depending
on the process used to generate or transmit the digital image.
For obvious reasons, it is desirable to remove the noise in a visually plausible manner.

\subsection{Anisotropic Diffusion}\label{sect:aniso}
In PDE-based diffusion filtering we typically process the image by solving a PDE of the form
\begin{subequations}
	\begin{align}\label{eqn:ch2_3}
		& \pd_t u = \Div \left( \vec{j} \right) && \text{in} \; \Omega, t > 0,  \\
		& u|_{t=0} = u^0, \\
		& \vec{j}^\T \vec{\nu} = 0 && \text{on} \; \bd \Omega, t > 0 \;.
	\end{align}
\end{subequations} 
Here, $u^0$ is the initially given image, $\Omega$ denotes the rectangular image domain, $\vec{\nu}$ is the exterior boundary normal, and
$\vec{j}$ is the flux vector. For isotropic filters the flux vector takes the form
\begin{align*}
	\vec{j} &= \kappa \grad u \;.
\end{align*}
If we choose $\kappa=1$ here we have Gaussian diffusion.
If we choose the Perona--Malik function $\kappa = g(|\grad u|^2)$ we would have
Perona--Malik\footnote{In their paper, Perona and Malik referred to their method as anisotropic although we consider it to be isotropic in the sense that the flux is always parallel to the gradient.}
edge-stopping diffusion.

Depending on the choice of $\kappa$, isotropic filters
can blur an edge (e.g., Gaussian diffusion), displace an edge,
or just stop the diffusion there (e.g., Perona--Malik).
However isotropic filters cannot smooth noise along an edge.
This is because the flux is parallel to the gradient, which is perpendicular to edges.

An anisotropic flux vector could take the form
\begin{align*}
	\vec{j} &= \kappa_1 \grad u + \kappa_2 \grad u^{\perp}
\end{align*}
where $\kappa_1$ and $\kappa_2$ are chosen based on the local structure of the image.
For example, near an edge, it would be desirable to prevent smoothing across the edge ($\kappa_1$ small) and instead smooth only along the edge ($\kappa_2$ large).

Weickert's approach \cite{w98} to implement this idea was to design a diffusion tensor
\begin{align}\label{eqn:Gtensor}
	\tens{G}[u] &= \kappa_1 \vec{\omega}^{\perp} \vec{\omega}^{\perp \T} + \kappa_2 \vec{\omega}  \vec{\omega}^\T, 
\end{align}
where $\vec{\omega} = \vec{\omega}[u]$ is a normalized ``edge vector'' that approximates the visual edge in the image $u$ in a small neighborhood about the point $x$.
The corresponding flux is given by
\begin{align}\label{eqn:Gtensorflux}
	\vec{j} &= \tens{G}[u] \grad u \; = \kappa_1 \left(\vec{\omega}^{\perp \T}\grad u \right) \; \vec{\omega}^{\perp}  + \kappa_2 \left(\vec{\omega}^\T \grad u \right) \; \vec{\omega}
\end{align}
and is directed by the visual edge $\vec{\omega}$ rather than the gradient direction.
The tensor $\tens{G}[u] \in \Real^{2 \times 2}$ is a nonlinear function of $u$ and varies spatially with $x$.
The diffusion is anisotropic since $\tens{G}$ is not simply a scalar multiple of the identity matrix.

The visual edge vector $\vec{\omega}$ needed to construct $\tens{G}$ is found by \emph{structure tensor analysis} \cite{w98}.
The structure tensor and the different choices for $\kappa_1,\kappa_2$ in \eqref{eqn:Gtensor}---which lead to edge- and coherence-enhancing diffusion---are discussed briefly in
the following.

\subsection{The Structure Tensor}\label{sec:structure_tensor}
The structure tensor provides information about orientations and coherent structures in an image \cite{w98,bwbm06,aach2006analysis}.
In its simplest form, the structure tensor is the matrix
\begin{align*}
	\tens{J}_{\sigma,0}[u] &:= \grad u_{\sigma} \grad u_{\sigma}^{\text{T}},
\end{align*}
where---in order to reduce the effect of noise or irrelevant small scale features---the image $u$ has been smoothed with a heat-kernel
\begin{align*}
	u_{\sigma} &:= K_{\sigma} \ast u, & K_\sigma(x) &= \frac{e^{\left(-\frac{x^2}{4\sigma} \right)}}{\sqrt{4\pi\sigma}}.
\end{align*}
In contrast to \cite{w98} we prefer the use of heat kernels $K_\sigma$ over Gaussians as this will extend in a straight-forward manner to surfaces (see Section~\ref{sect:aniso_surf}), and has the same effect in practice.
$\tens{J}_{\sigma,0}[u]$ is the initial form of the structure tensor;
coherent structures over a bigger neighborhood around the point $x$ are often accounted for by averaging tensor-valued orientation information given in $\tens{J}_{\sigma,0}$.
This is done by component-wise convolution of $\tens{J}_{\sigma,0}$ with a second heat-kernel $K_\rho$
to define the structure tensor
\begin{align}
	\tens{J}_{\sigma,\rho}[u] &:= K_{\rho} \ast \tens{J}_{\sigma,0}[u] = K_{\rho} \ast \left( \grad u_{\sigma} \grad u_{\sigma}^\T \right).
\end{align}
By construction, $\tens{J}_{\sigma,\rho}[u]$ is symmetric positive semidefinite, and has therefor a spectral decomposition
\begin{align}\label{eqn:Jspectral}
	\tens{J}_{\sigma,\rho}[u] &= \mu_1 \vec{\omega}^{\perp} \vec{\omega}^{\perp \T} + \mu_2 \vec{\omega}  \vec{\omega}^\T, \quad \mu_1 \geq \mu_2 \geq 0.
\end{align}
This decomposition is unique if $\mu_1 \neq \mu_2$ and yields the orientations parallel $\vec{\omega}$ and perpendicular $\vec{\omega}^\perp$ to the visual edge.
(Note that if $\rho=0$, one obtains $\mu_2 = 0$ and $\vec{\omega} || \grad u_{\sigma}^{\perp}$.)
The eigenvalues give information about the image structure near $x$, specifically (cf. \cite{w98}) 
\begin{subequations} 
	\begin{align}
		&\mu_1 = \mu_2 && \implies x \; \text{in homogeneous region}, \\
		&\mu_1 \gg \mu_2 = 0 && \implies x \; \text{near straight edge}, \\
		&\mu_1 \ge \mu_2 \gg 0 && \implies x \; \text{near corner}.
	\end{align}
\end{subequations} 
The so-called \emph{coherence} $c := \mu_1 - \mu_2 \geq 0$ measures how pronounced the edge is.
It also gives a measure for the local contrast and tells us if the decomposition \eqref{eqn:Jspectral} is well- or ill-conditioned.

Now, we consider the choices of $\kappa_1$, $\kappa_2$ in \eqref{eqn:Gtensor} to obtain the diffusion tensor $\tens{G}[u]$.
The eigenvectors $\vec{\omega}^{\perp}$, $\vec{\omega}$ of $\tens{G}[u]$ are exactly those of the structure tensor $\tens{J}_{\sigma,\rho}[u]$ from \eqref{eqn:Jspectral},
while we replace the eigenvalues with $\kappa_1$, $\kappa_2$, depending on the coherence $c$.

\subsection{Edge-Enhancing Diffusion}
In edge-enhancing diffusion, the goal is to smooth noise while keeping or enhancing edges.
According to Weickert \cite{w98}, one works with the initial structure tensor $\tens{J}_{\sigma,0}[u]$.
The appropriate choice of $\kappa_1$ and $\kappa_2$  is based on the following observation:
in homogeneous regions $|\grad u_{\sigma}|^2 =  \mu_1 \approx \mu_2 = 0$, whereas near an edge $\mu_1 \gg \mu_2 =0$.
Thus, $\kappa_1$ (the diffusivity across edges) is chosen to be a decreasing function of $\mu_1$, 
while $\kappa_2$ (the diffusivity along edges) is set equal to one:
\begin{align}\label{eqn:ee_ews}
	\kappa_1 = g(\mu_1) \quad \text{and} \quad \kappa_2 = 1.
\end{align}
Here $g$ is a scalar-valued function with $g(0) = 1$ and $\lim_{s \to \infty} g(s) = 0$.
We choose the Perona--Malik diffusivity function
\begin{align}\label{eqn:g}
	g(s^2) = \frac{1}{1 + \frac{s^2}{\lambda^2}},
\end{align}
where $\lambda \in \Real$ in is tunable parameter that determines the filter's sensitivity to edges.
For a different choice of $g$ see \cite{w98}.
In regions of very little contrast, we have $\kappa_1 \approx 1$ and the diffusion will behave locally like Gaussian diffusion since $\tens{G}[u] \approx I$.
Near an edge, we have $\kappa_1 \approx 0$ which implies $\tens{G}[u] \approx \vec{\omega} \vec{\omega}^{\text{T}}$, and diffusion occurs mainly along the edge.

\begin{figure}
        \centering
        \includegraphics[width=.48\columnwidth]{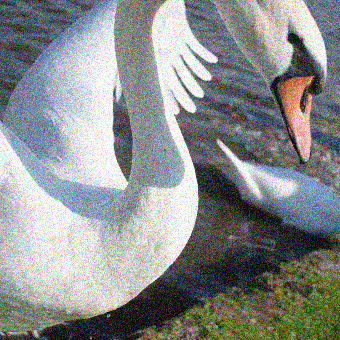}\;\;%
        \includegraphics[width=.48\columnwidth]{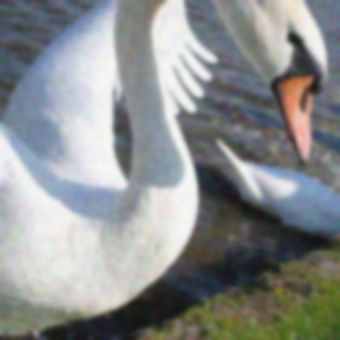}\vspace{1ex}\\
        \includegraphics[width=.48\columnwidth]{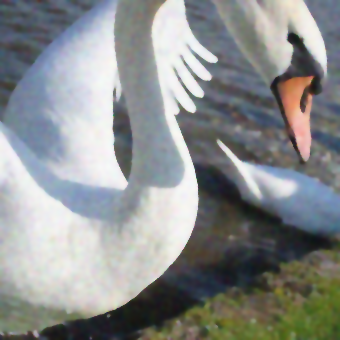}\;\;%
        \includegraphics[width=.48\columnwidth]{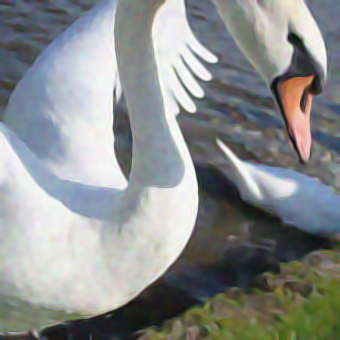}%
\caption{The effect of different diffusion filters.
\textit{Top Left:} Noisy image (photograph by courtesy of Harry Biddle).
\textit{Top Right:} Gaussian diffusion.
\textit{Bottom Left:} Perona--Malik diffusion; edge-sensitivity parameter $\lambda_\rel = 0.3$.
\textit{Bottom Right:} Edge-enhancing diffusion; $\sigma = 0.75$, $\rho = 2.5$, $\lambda_\rel = 10^{-2}$.
Stop time in all cases $T=2.5$.}
\label{fig:swan_comp}
\end{figure}

The edge-enhancing set-up of $\tens{G}[u]$ can also be based on $\tens{J}_{\sigma,\rho}[u]$.
But with $\rho>0$, $\kappa_1 = g(c)$ is a function of the coherence $c := \mu_1 - \mu_2$. Because the decomposition \eqref{eqn:Jspectral} is ill-conditioned if $c \approx 0$,
we force $\tens{G}[u] \approx \tens{I}$ in that case.
This is also compatible with the case $\rho=0$, where $c = \mu_1$ because $\mu_2=0$.
An example of edge-enhancing diffusion filtering is shown in Figure~\ref{fig:swan_comp}.

\subsection{Coherence-Enhancing Diffusion}
The second type of anisotropic diffusion which we consider is coherence-enhancing diffusion which
aims to preferentially smooth in the direction of largest coherence. 
Weickert \cite{w98} suggests the following choice for $\kappa_1$ and $\kappa_2$ in \eqref{eqn:Gtensor}:
\begin{align}\label{eqn:ce_ews}
	\kappa_1 &= \alpha, & \kappa_2 &= \alpha + (1 - \alpha) \exp \left( \dfrac{-B^2}{(\mu_1 - \mu_2)^2} \right), 
\end{align}
where $0 \le \alpha \ll 1$ and $B>0$ are constants. $\kappa_2$ is a function of the coherence $c=\mu_1-\mu_2$ where
$\mu_1$ and $\mu_2$ are the eigenvalues of $\tens{J}_{\sigma,\rho}[u]$. 
Diffusion along the edge ($\kappa_2$ large) is performed when the coherence is large compared to $B$.
Where the coherence is small, i.e., $\kappa_2 \approx \kappa_1 = \alpha$, only a small amount of Gaussian diffusion is performed.

\begin{figure}
	\centering
		\includegraphics[width=0.91\columnwidth]{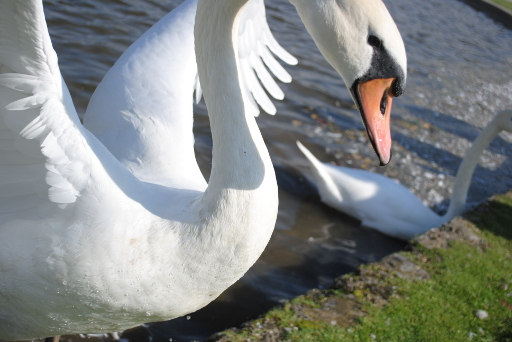}\vspace{1ex}\\
		\includegraphics[width=0.91\columnwidth]{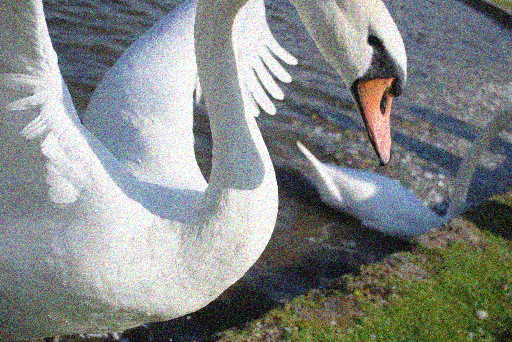}\vspace{1ex}\\
		\includegraphics[width=0.91\columnwidth]{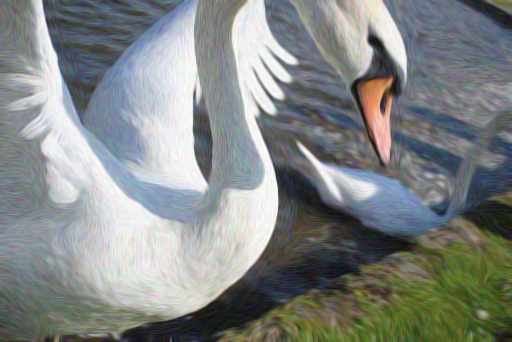}\vspace{1ex}\\
		\includegraphics[width=0.91\columnwidth]{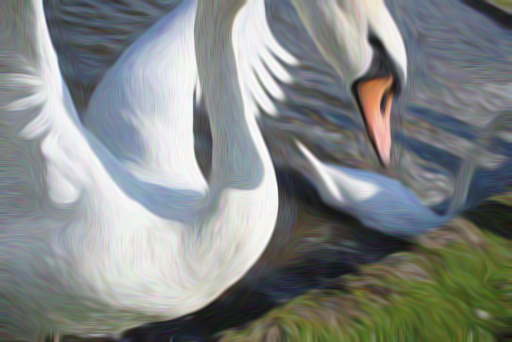}
	\caption{Stylized texture created by coherence-enhancing diffusion. From top to bottom:
				\textit{1.} Original photograph.
				\textit{2.} Noise added for texture creation.
				\textit{3.} Stop time $T = 2.5$. 
				\textit{4.} Stop time $T = 7.5$. 
				The parameters are $\sigma = 0.75$, $\rho = 25$, $\alpha = 10^{-3}$, and $B_\rel = 10^{-10}$.}
	\label{fig:swan_texture}
\end{figure}

One application of coherence-enhancing diffusion lies in its ability to create stylized images (cf. \cite{w98})
which might be desired for use in comic books or video games.
Coherence-enhancing diffusion creates a natural looking texture because of smoothing according to the inherent structure of the image.
As the stop time increases, the image becomes less photorealistic and more stylized; this is shown in Figure~\ref{fig:swan_texture}.

\subsection{Color Images}
If the given image $u^0$ is a color image the diffusion is applied separately to each channel $u_k$
\begin{align}\label{eqn:colorIm}
	\pd_t u_k &= \Div \left( \tens{G}[u] \grad u_k \right) \;.
\end{align}
In order to minimize spurious colors, we use the same diffusion tensor $\tens{G}[u]$ for all channels which is found (as in the scalar case) from a common structure tensor.
We compute the common structure tensor following \cite{weickert1999coherence} as
\begin{align}\label{eqn:CommonST}
	\tens{J}_{\sigma,\rho}[u] &:= \frac{1}{n} \sum\limits_{k=1}^{n} \tens{J}_{\sigma,\rho}[u_k],
\end{align}
i.e., the arithmetic mean of all channel structure tensors $\tens{J}_{\sigma,\rho}[u_k]$.
Depending on the color space some channels might be more important to structure information than others.
In this case a general weighted average can be used which puts larger weights on the important channels.

\section{Surface-Intrinsic Anisotropic Diffusion}\label{sect:aniso_surf}
In this section we formulate Weickert's diffusion models for images defined on surfaces.
We consider two-dimen\-sional surfaces $\surf$ embedded in $\Real^3$ which are smooth
and orientable with a normal field $\vec{n}$.
Images are now functions $u:\surf \to \Real$ defined on $\surf$.
The diffusion PDE to process the image takes the form
\begin{subequations}
	\begin{align}
		& \pd_t u = \sDiv \left( \tens{G}[u] \sgrad u \right) && \text{in} \; \surf, t > 0,  \label{eqn:GSdiff} \\
		& u|_{t=0} = u^0, \\
		& (\tens{G}[u] \sgrad u) ^\T \vec{\nu} = 0 && \text{on} \; \bd \surf, t > 0 \;. \label{eqn:GSdiffBC}
	\end{align}
\end{subequations}
The boundary condition \eqref{eqn:GSdiffBC} is relevant only in cases where the surface is open.
In this case $\vec{\nu}$ is an exterior co-normal, i.e., it is perpendicular to the boundary curve but tangent to the surface $\vec{\nu}^\T \vec{n}=0$.

The differential operators in PDE \eqref{eqn:GSdiff} are intrinsic to the surface and are related to the standard operators as follows (see e.g., \cite{Maerz/Macdonald:cpfunctions})
\begin{subequations}\label{eqn:sgraddiv}
	\begin{align}
		\sgrad u &= \grad u - (\vec{n}^\T \grad u) \vec{n}, \label{eqn:sgrad} \\
		\sDiv (\vec{j}) &= \Div(\vec{j}) - \vec{n}^\T D\vec{j} \; \vec{n}, \label{eqn:sDiv}
	\end{align}
\end{subequations}
where $D\vec{j}$ is the Jacobian matrix of the flux $\vec{j}$. 
Note that, with this form, the gradient $\sgrad u$ and the flux $\vec{j}$ are both three-component vectors 
while $\tens{G}[u]$ is a symmetric ${3 \times 3}$-matrix.  

While not strictly required, it may be desirable that the PDE \eqref{eqn:GSdiff} represent a surface conservation law (for example, to guarantee that at least the continuous problem conserves overall gray levels).
This means that the flux vector $\vec{j} = \tens{G}[u] \sgrad u \in \Real^3$ should always be tangent to the surface $\surf$, i.e., $\vec{j}^\T \vec{n}=0$.
This is because the divergence theorem holds only for the tangential part $\vec{j}_{\tan}$:
\begin{align}\label{eqn:divtheo}
	\int\limits_{\bd \Omega} \vec{j}_{\tan}^\T \vec{\nu} d\len &= \int\limits_{\Omega} \sDiv (\vec{j}_{\tan}) d\area,
\end{align}
where $\Omega \subset \surf$, $d\len$ and $d\area$ are the infinitesimal length and area measures.
Since the co-normal is tangent to the surface we have $\vec{j}^\T \vec{\nu} = \vec{j}_{\tan}^\T \vec{\nu}$,
but \eqref{eqn:divtheo} does not hold when $\vec{j}_{\tan}$ is replaced with $\vec{j}$ (unless they are equal).
Thus only a tangential flux $\vec{j} = \vec{j}_{\tan}$ will turn PDE \eqref{eqn:GSdiff} into a surface conservation law.

The design of a diffusion tensor $\tens{G}[u]$ which obeys Weickert's model 
and maps vectors to the tangent space of $\surf$ to give a tangential flux is the topic of the next section.  

\subsection{The Visual Edge Vector on a Surface}\label{sec:edge_vector}
We begin with finding visual-edge vectors for functions $u$ defined on a surface.
We leverage the approach of \cite{aach2006analysis} to find a surface intrinsic form of the structure tensor.
Let $\vec{\omega}$ again denote the edge vector and let $\Omega(x) \subset \surf$ denote a neighborhood of the point $x \in \surf$, then we require
\begin{subequations}
	\begin{align}
		(\vec{\omega}^\T \sgrad u(y))^2 &\approx 0 \quad \forall \; y \in \Omega(x), \label{eqn:almostperp}\\
		|\vec{\omega}|^2 &= 1, \label{eqn:nontriv}\\
		\vec{\omega}^\T \vec{n}(x) &= 0 .\label{eqn:tangent}
	\end{align}
\end{subequations}
Requirement \eqref{eqn:almostperp} says that we want a vector $\vec{\omega}$ which is almost perpendicular to gradients in a neighborhood of $x$.
Condition \eqref{eqn:nontriv} excludes trivial solutions $\vec{\omega} = 0$. Finally, $\vec{\omega}$ must be tangent to the surface at $x$ which is expressed in \eqref{eqn:tangent}.

As in \cite{aach2006analysis}, we formulate \eqref{eqn:almostperp} as least squares problem.
Let $K_\rho(x,y)$ denote the surface heat kernel on $\surf$ with $\rho$ proportional to
the diameter of $\Omega(x)$, then we have
\begin{align}\label{eqn:LSform}
	\vec{\omega} &= \arg \min\limits_{\vec{v}} \sInt K_\rho(x,y) (\vec{v}^\T \sgrad u(y))^2 \; d\area(y),
\end{align}
where the minimization happens over all vectors $\vec{v}$ satisfying \eqref{eqn:nontriv} and \eqref{eqn:tangent}.
With the $3 \times 3$ structure tensor
\begin{align}
	\tens{J}_{0,\rho}[u](x) &:= \sInt K_\rho(x,y) \sgrad u(y) \sgrad u(y)^\T \; d\area(y)
\end{align}
we can rewrite \eqref{eqn:LSform} as
\begin{align}
	\vec{\omega} &= \arg \min\limits_{\vec{v}} \vec{v}^\T \; \tens{J}_{0,\rho}[u](x) \; \vec{v}.
\end{align}
If we consider $u_\sigma$, i.e., $u$ smoothed with the surface heat kernel $K_\sigma(x,y)$, in place of $u$ in \eqref{eqn:almostperp} we will get
\begin{align}\label{eqn:Jdef}
	\tens{J}_{\sigma,\rho}[u](x) &:= \sInt K_\rho(x,y) \sgrad u_\sigma(y) \sgrad u_\sigma(y)^\T \, d\area(y)
\end{align}  
as the structure tensor. 
The new problem is now
\begin{subequations}
	\begin{align}
		\vec{\omega} &= \arg \min\limits_{\vec{v} \in \Real^3} \vec{v}^\T \; \tens{J}_{\sigma,\rho}[u](x) \; \vec{v}, \label{eqn:LSST}\\
		|\vec{v}|^2 &= 1, \label{eqn:nontriv2}\\
		\vec{v}^\T \vec{n}(x) &= 0, \label{eqn:tangent2}
	\end{align}
\end{subequations}
and can be solved with Lagrange multipliers which leads to the following condition
\begin{align}\label{eqn:Lag}
	\tens{J} \; \vec{v} - \lambda \vec{v} - \mu \vec{n} &= 0
\end{align}  
subject to the side conditions \eqref{eqn:nontriv2} and \eqref{eqn:tangent2}.
Here we use $\tens{J}$ and $\vec{n}$ as short hands for $\tens{J}_{\sigma,\rho}[u](x)$ and $\vec{n}(x)$.

In order to reduce the system, we use an orthonormal basis of the tangent space $\Tang{x} = \lingen \{\vec{q}_1(x),\vec{q}_2(x)\}$.
For $\tilde{\vec{v}} \in \Real^2$ we set
\begin{align}
	\vec{v} &= \tens{Q} \tilde{\vec{v}} \quad \text{with} \quad \tens{Q} := \left[ \vec{q}_1 | \vec{q}_2 \right] \in \Real^{3 \times 2},
\end{align}
then \eqref{eqn:tangent2} is satisfied and condition \eqref{eqn:nontriv2} implies $|\tilde{\vec{v}}|^2=1$.
Moreover, \eqref{eqn:Lag} turns into
\begin{align}\label{eqn:Lag2}
	\tens{J} \; \tens{Q} \tilde{\vec{v}} - \lambda \tens{Q} \tilde{\vec{v}} - \mu \vec{n} &= 0.
\end{align}  
Finally, left multiplication with $\tens{Q}^\T$ reduces \eqref{eqn:Lag2} to a $2 \times 2$ eigenvalue problem
\begin{align}\label{eqn:Lag3}
	\tens{Q}^\T \tens{J} \; \tens{Q} \tilde{\vec{v}} - \lambda \tilde{\vec{v}} &= 0.
\end{align}
The $2 \times 2$ tensor $\tilde{\tens{J}}$ given by the tensor contraction
\begin{align}
	\tilde{\tens{J}} &:= \tens{Q}^\T \tens{J} \; \tens{Q}
\end{align}
can be seen as the surface-intrinsic structure tensor. Now, the spectral decomposition of $\tilde{\tens{J}}$
\begin{align}
	\tilde{\tens{J}} &= \mu_1 \tilde{\vec{\omega}}^\perp  \tilde{\vec{\omega}}^{\perp \T} + \mu_2  \tilde{\vec{\omega}} \tilde{\vec{\omega}}^\T, \quad \mu_1 > \mu_2
\end{align}
yields the desired solution: the eigenvector $\tilde{\vec{\omega}} \in \Real^2$ with respect to the minimal eigenvalue $\mu_2$ represents the visual edge in tangent space coordinates,
hence in embedding space coordinates we have
\begin{align}
	\vec{\omega} &= \tens{Q} \tilde{\vec{\omega}} \quad \text{and} \quad \vec{\omega}^\perp = \tens{Q} \tilde{\vec{\omega}}^\perp.
\end{align}

The surface-intrinsic $2 \times 2$ diffusion tensor $\tilde{\tens{G}}$ is constructed by
\begin{align}
	\tilde{\tens{G}} &:= \kappa_1 \tilde{\vec{\omega}}^\perp  \tilde{\vec{\omega}}^{\perp \T} + \kappa_2  \tilde{\vec{\omega}} \tilde{\vec{\omega}}^\T,
\end{align}
where we choose $\kappa_1,\kappa_2$ as discussed in Section~\ref{sect:aniso} to get either edge- or coherence-enhancing diffusion.
Finally, the $3 \times 3$ version $\tens{G}$, embedded in $\Real^3$ to be used with PDE \eqref{eqn:GSdiff}, is obtained by reverting the tensor contraction
\begin{align}
	\tens{G} = \tens{Q} \tilde{\tens{G}} \; \tens{Q}^\T.
\end{align}
Note that the flux $\vec{j} = \tens{G} \sgrad u$ will automatically be tangential to $\surf$ since the columns of $\tens{Q}$ are a basis of the tangent space.

In \eqref{eqn:Jdef} we based the definition  $\tens{J}_{\sigma,\rho}[u](x)$ on integration with respect to surface heat kernels.
Instead, in practice, we solve the surface-intrinsic Gaussian diffusion equation
\begin{subequations}
	\begin{align}
		& \pd_\tau w = \slap w && \text{in} \; \surf, \tau > 0,  \label{eqn:heat} \\
		& w|_{\tau=0} = w^0, \\
		& \sgrad w^\T \vec{\nu} = 0 && \text{on} \; \bd \surf, \tau > 0 \;. \label{eqn:NeuBC}
	\end{align}
\end{subequations}
where $\slap$ is the Laplace-Beltrami operator. 
For the pre-smoothing step we set $w^0 = u$ and solve until $\tau = \sigma$, for the post-smoothing step
we initialize $w^0$ with the $ij$-component of $\tens{J}_{\sigma,0}[u]$ and solve until $\tau=\rho$.
Since $\tens{J}_{\sigma,\rho}[u]$ is symmetric, this means six further heat solves. 

Algorithm~\ref{alg:G3D} summarizes all the steps to get the diffusion tensor $\tens{G}$.
\begin{algorithm}[H]
	\caption{Construction of $\tens{G}[u]$}\label{alg:G3D}
	\begin{algorithmic}[1]
		\State Obtain $\tens{Q}$ as described in Section~\ref{sect:Tang}.
		\Comment{(Step 1 can be done once prior to evolving PDE \eqref{eqn:GSdiff})}
		\State Find $u_{\sigma}$ with a heat solve until $\sigma$.
		\State Calculate $\sgrad u_{\sigma}$ and initialize $\tens{J}_{\sigma,0}[u] = \sgrad u_{\sigma} \sgrad u_{\sigma}^{\text{T}}$
		\State Find $\tens{J}_{\sigma,\rho}[u]$  with component-wise heat solves until $\rho$
		\State Calculate $\tilde{\tens{J}}_{\rho,\sigma}[u] = \tens{Q}^\T \tens{J}_{\sigma,\rho}[u] \tens{Q}$ 
		\State Find the spectral decomposition of $\tilde{\tens{J}}_{\rho,\sigma}[u]$
		\State Set up $\tilde{\tens{G}}[u]$ according to edge-enhancing or coherence-enhancing diffusion
		\State Calculate $\tens{G}[u] = \tens{Q} \tilde{\tens{G}}[u] \; \tens{Q}^\T$
	\end{algorithmic}
\end{algorithm}

\subsection{Finding the Tangent Space Basis}\label{sect:Tang}
We can offer two point-wise approaches to finding an orthonormal basis of the tangent space $\Tang{x}$ giving
the columns of $\tens{Q}(x) = \left[\vec{q}_1(x) | \vec{q}_2(x) \right]$.
If $\vec{n}$ is given as part of the problem description, then for each point $x$ on the surface, we perform
a $QR$-decomposition (here a single Householder reflection) of the vector $\vec{n}(x)$:
\begin{align}
	\vec{n}(x) &= \left[ \pm \vec{n}(x) | \vec{q}_1(x) | \vec{q}_2(x) \right] \; (\pm 1,0,0)^\T.
\end{align}
The second and third column are the desired orthonormal tangent vectors.

In Section~\ref{sec:cpm} we will introduce the closest point method which uses a closest point function $\cp$ to represent the surface.
In \cite{Maerz/Macdonald:cpfunctions} we showed that the Jacobian $D\cp(x)$ when evaluated at surface points $x$ yields the orthogonal projection matrix $\tens{P}(x)$ 
which projects onto the tangent space $\Tang{x}$, i.e.,
\begin{align}\label{eqn:DcpP}
	D\cp(x) &= \tens{P}(x) = \tens{I} - \vec{n}(x) \vec{n}(x)^T \qquad x \in \surf.
\end{align}
Because of \eqref{eqn:DcpP} we can find $\vec{q}_1(x)$ and $\vec{q}_2(x)$ from the spectral decomposition of $D\cp(x)$:
\begin{align}\label{eqn:Dcp}
	D\cp(x) &= \left[ \vec{n} | \vec{q}_1 | \vec{q}_2 \right](x) \,
	\left[
	\begin{matrix}
		0 & 0 & 0\\
		0 & 1 & 0\\
		0 & 0 & 1
	\end{matrix}
	\right] \,
	\left[ \vec{n} | \vec{q}_1 | \vec{q}_2 \right](x)^T.
\end{align}

Now, we have all ingredients of the surface-intrinsic diffusion model.
It involves several linear and non-linear in-surface diffusion steps.
In the following section we explain how we deal with these equations in order to solve them numerically.

\section{Solving In-Surface Diffusion Equations with the Closest Point Method} \label{sec:cpm}
The closest point method, introduced in \cite{Ruuth/Merriman:jcp08:CPM}, is an embedding technique for solving PDEs posed on embedded surfaces.
The central idea is to extend functions and differential operators to the surrounding space (here $\Real^3$) and to solve an embedding equation which is
a 3D-analog of the original surface equation. This approach is appealing since the extended versions of the operators $\sgrad$, $\sDiv$, and $\slap$
in the closest point framework turn out to be $\grad$, $\Div$, and $\lap$.
Hence, embedding equations are accessible to standard finite difference techniques and existing algorithms (and even codes) can be reused.

\subsection{Surface Representation}
The closest point method utilizes the \textit{closest point representation} of a surface which is given in terms of the closest point function
\begin{align}\label{eqn:cp_function}
	\cp(x) = \arg \min\limits_{\hat{x} \in \surf} |x - \hat{x}|.
\end{align}
For a point $x$ in the embedding space, $\cp(x)$ is the point on the surface $\surf$ which is closest in Euclidean distance to $x$.
This function is well-defined in a tubular neighborhood or narrow band $\band$ of the surface and is as smooth as the underlying surface \cite{Maerz/Macdonald:cpfunctions}.

\begin{figure}
	\centering
	\definecolor{darkgreen}{rgb}{0,.6,0}
	\definecolor{lightcyan}{rgb}{0,.75,.75}
	\begin{tikzpicture}[scale=0.5]

	\draw[lightcyan,very thick] (0.1,0.3) arc (160:97:10 and 7.3);

	\filldraw[white] (0,0) circle (.1pt);

	\def\gdsz{2.5pt}
	\foreach \i/\x/\y\t in {
		1/0/1/156,
		2/2/1/146,
		3/4/1/133,
		4/6/1/117.4,
		5/0/3/146.6,
		6/2/3/137,
		7/4/3/125.9,
		8/6/3/113,
		8/8/3/100.5,
		10/0/5/139,
		11/2/5/130.5,
		12/4/5/120.6,
		13/6/5/110,
		14/8/5/98.64,
		15/2/7/126.1,
		16/4/7/117,
		17/6/7/107,
		18/8/7/97.5
	}{
		\path (0.1,0.3) arc (160:\t:10 and 7.3) node (v\i) {}; 
		\draw[latex-] (v\i.center) -- (\x,\y);
		\filldraw[fill=white] (\x,\y) circle (\gdsz);
	};

	\end{tikzpicture}
	\caption{Closest point representation of a curve embedded in $\Real^2$.
				For each point $x$ (black circles) in the embedding space, the point (arrow tip) on the cyan curve, which is closest in Euclidean distance to $x$, is stored.}
				\label{fig:cp_grid}
\end{figure}
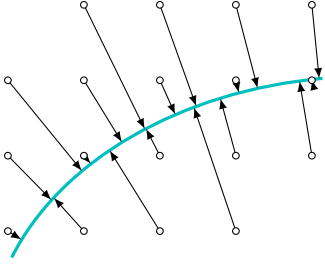

The closest point function can be derived analytically for many common surfaces.
When a parameterization or triangulation of the surface is known, the closest point function can be computed numerically \cite{Ruuth/Merriman:jcp08:CPM}.
Figure \ref{fig:cp_grid} illustrates the closest point representation of a curve embedded in $\Real^2$. 

\subsection{Extending Functions and Differential Operators}
Using the closest point representation, we can extend values of a surface function $u:\surf \to \Real$ 
into the surrounding band $\band$ by defining $\bar{u} : \band \to \Real$ as
\begin{align}\label{eqn:cpext}
  \bar{u}(x) &:= u(\cp(x)).
\end{align}
Notably, $\bar{u}$ will be constant in the direction normal to the surface and this property is key to the closest point method:
it implies that an application of a Cartesian differential operator to $\bar{u}$
is equivalent to applying the corresponding intrinsic surface differential operator to $u$.
We state this as principles below; these mathematical principles were established in \cite{Ruuth/Merriman:jcp08:CPM} and proven in \cite{Maerz/Macdonald:cpfunctions}.
\begin{principle} \label{thm:gradient_principle}
	\emph{(Gradient Principle):} Let $\surf$ be a surface embedded in $\Real^n$ and let $u$ be a function, defined on $\Real^n$,
	that is constant along directions normal to the surface, then
	\begin{align}
		\grad u(x) &= \sgrad u(x) \qquad \forall \; x \in \surf. 
	\end{align}
\end{principle}
\begin{principle} \label{thm:divergence_principle}
	\emph{(Divergence Principle):} Let $\surf$ be a surface embedded in $\Real^n$. If $\vec{j}$ is a vector field on $\Real^n$ that is
	tangent to $\surf$ and tangent to all surfaces displaced a fixed Euclidean distance from $\surf$, then
	\begin{align}
		\Div \vec{j}(x) &= \sDiv \vec{j}(x) \qquad \forall \; x \in \surf.
	\end{align}
\end{principle}
Direct consequences of these principles are
\begin{align}
	\sgrad u(x) &= \grad \bar{u}(x), \label{4D_first} \\
	\sDiv \vec{j}(x) &= \Div \bar{\vec{j}}(x), \qquad \forall \; x \in \surf, \label{4D_second}
\end{align}
where $\bar{\vec{j}}$ is the closest point extension of a tangential flux $\vec{j}$.
Moreover, since $\grad \bar{u}$ is tangential to level-surfaces of the Euclidean distance-to-$\surf$ map \cite{Ruuth/Merriman:jcp08:CPM,Maerz/Macdonald:cpfunctions}, 
combining Principles~\ref{thm:gradient_principle} and \ref{thm:divergence_principle} yields
\begin{align}
	\slap u(x) &= \lap \bar{u}(x), \label{4D_third} \\
	\sDiv \left( g \sgrad u \right)(x) &= \Div \left( \bar{g} \grad \bar{u} \right)(x), \qquad \forall \; x \in \surf, \label{4D_third2}
\end{align}
where $g:\surf \to \Real$ is a scalar diffusivity and $\bar{g}$ its closest point extension.

Additionally, if $\tens{G}$ is a diffusion tensor which maps to the tangent space, i.e., the corresponding flux $\vec{j} = \tens{G} \sgrad u$ is tangential, then
Principles \ref{thm:gradient_principle} and \ref{thm:divergence_principle} also imply that
\begin{align}
	\sDiv \left( \tens{G} \sgrad u \right)(x) &= \Div \left( \bar{\tens{G}} \grad \bar{u} \right)(x), \qquad \forall \; x \in \surf, \label{4D_fourth}
\end{align}
where $\bar{\tens{G}}$ is the closest point extension of $\tens{G}$.

Finally, we note that a closest point extension $\bar{u}$ is characterized \cite{vonGlehn:mol} by
\begin{align}
	\bar{u} &= \bar{u} \circ \cp,
\end{align}
i.e., it is the closest point extension of itself.

\subsection{Gaussian Diffusion}
We start with the Gaussian diffusion equation on a closed surface $\surf$ in order to demonstrate the embedding idea:
\begin{subequations}
	\begin{align}
		& \pd_t u = \slap u && x \in \surf, \; t > 0, \label{4B_first} \\
		& u|_{t=0} = u^0. \label{4B_second}
	\end{align}
\end{subequations}  
Using \eqref{4D_third} we obtain the following embedding problem
\begin{subequations}  \label{4C}
	\begin{align}
		& \pd_t v = \lap v && x \in \band, \; t > 0, \label{4C_first} \\
		& v|_{t=0} = u^0 \circ \cp \label{4C_second} \\
		& v = v \circ \cp && x \in \band, \; t > 0, \label{4C_third} \;.
	\end{align}
\end{subequations}  
Here \eqref{4C_second} says that we start the process with a closest point extension of the initial data $u^0$, while
condition \eqref{4C_third} guarantees that $v$ stays a closest point extension for all times and hence we can rely on the principles
which give us \eqref{4C_first} as the 3D-analog of \eqref{4B_first}. 

In order to cope with condition \eqref{4C_third} Ruuth \& Merriman \cite{Ruuth/Merriman:jcp08:CPM} suggested the following semi-discrete (in time) iteration:
after initialization $v^0 = u^0 \circ \cp$, alternate between two steps
\begin{subequations}
	\begin{align}
		&1.\text{ Evolve} & \quad w &= v^n + \dt \lap v^n, \\
		&2.\text{ Extend} & \quad v^{n+1} &= w \circ \cp,
	\end{align}
\end{subequations}
where $\dt$ is the time-step size. 
Here step 1 evolves \eqref{4C_first} of the embedding problem, while step 2 reconstructs the surface function or rather its closest point extension
to make sure that \eqref{4C_third} is satisfied at time $t_{n+1}$. 

A fully discrete version of the Ruuth \& Merriman iteration needs a discretization of the spatial operators in step 1 
and an interpolation scheme in step 2. Using a uniform Cartesian grid in $\Real^3$
the Laplacian in our example can be discretized with the standard 7 point finite difference stencil
and be implemented as a matrix $L$ acting on the 1D array $v^{n+1}$.
The interpolation scheme in step 2 is necessary since the data $w$ is given only on grid points and $\cp(x)$ is hardly ever a grid point (even though $x$ is).
In order to get around that we interpolate the data $w$ with tri-cubic interpolation and extend the interpolant $W$ rather than $w$:
\begin{align}
	v^{n+1}_i &= W \circ \cp(x_i)
\end{align}
where $x_i$ is the grid point corresponding to the $i$-th component of the array $v^{n+1}$. Since tri-cubic interpolation is linear in the data $w$, this
operation can also be implemented as an extension matrix $E$ acting on the 1D array $w$. The fully discrete Ruuth \& Merriman iteration reads then
\begin{subequations}
	\begin{align}
		&1.\text{ Evolve} & \quad w &= v^n + \dt L v^n, \\
		&2.\text{ Extend} & \quad v^{n+1} &= E w.
	\end{align}
\end{subequations}

The computation is performed in a \emph{computational band} for two reasons:
 firstly, the closest point function is defined in the band $\band$
and hence we need sufficient resolution within $\band$ to resolve the
geometry of the surface.
Secondly, the code can be made more efficient by working on a narrow band surrounding the surface $\surf$.
A nice feature of the Ruuth \& Merriman iteration is that no artificial boundary conditions on the boundary of the band need to be imposed \cite{Ruuth/Merriman:jcp08:CPM}.
This has to do with the extension step: values at grid points are overwritten at each time step with the value at their closest points.
Note also, that no artificial boundary conditions are imposed in the embedding problem \eqref{4C}.

For the sake of efficiency optimization of the width of the computational band is reasonable. 
The bandwidth depends on the degree of the interpolant and on the finite difference stencil used. 
Suppose we use an interpolant of degree $p$ and are working in $d$-dimensions, this requires $(p+1)^d$ points around an interpolation point $\cp(x_i)$ \cite{Ruuth/Merriman:jcp08:CPM}.
Furthermore, each of the points in the interpolation stencil must be advanced in time with a finite difference stencil. 
As a rule of thumb the diameter of the convolution of the interpolation and finite difference stencil gives a good value for the bandwidth.
More details on finding the optimal band are given in \cite[Appendix A]{mr10}.

\subsection{Anisotropic Diffusion}\label{sec:numerical_scheme}
In the case of anisotropic diffusion the closest point embedding idea yields the following embedding problem
\begin{subequations}
	\begin{align}
		& \pd_t v = \Div \left( \tens{G}[v] \circ \cp \grad v \right) && x \in \band, \; t > 0,  \\
		& v|_{t=0} = u^0 \circ \cp  \\
		& v = v \circ \cp && x \in \band, \; t > 0, 
	\end{align}
\end{subequations}  
which is again a 3D-analog of the surface PDE.
Analogously to Gaussian diffusion, we obtain a semi-discrete Ruuth \& Merriman iteration, i.e.,
after initialization $v^0 = u^0 \circ \cp$, we alternate between two steps
\begin{subequations}
	\begin{align}
		&1.\text{ Evolve} & \; w &= v^n + \dt \Div \left( \tens{G}[v^n] \circ \cp \grad v^n \right), \label{eqn:evolve_aniso}\\
		&2.\text{ Extend} & \; v^{n+1} &= w \circ \cp.
	\end{align}
\end{subequations}
From here we obtain the discrete iteration:
using Algorithm~\ref{alg:G3D} on the state $v^n$ (where all surface-intrinsic differentials and heat solves are dealt with the closest point method)
yields the diffusion tensor $\tens{G}[v^n] \circ \cp$. The extension (step 2) is realized as explained in the previous section.
Finally, we discretize the diffusion operator in \eqref{eqn:evolve_aniso} with the following formally second-order accurate scheme:
\begin{equation}\label{eqn:scheme_aniso}
\begin{aligned}
        \Div &\left( \tens{G}[v^n] \circ \cp \grad v^n \right) \approx \\
	& \quad\,
		D_x^- \left(A_x^+ \tens{G}_{11}  \ds  D_x^+ v^n \right)
	+ D_x^c \left(\tens{G}_{12}  \ds  D_y^c v^n \right) \\
	& + D_x^c \left(\tens{G}_{13}  \ds  D_z^c v^n \right)  
	+ D_y^c \left(\tens{G}_{12}  \ds  D_x^c v^n \right) \\
	& + D_y^- \left(A_y^+ \tens{G}_{22}  \ds  D_y^+ v^n \right)
	+ D_y^c \left(\tens{G}_{23}  \ds  D_z^c v^n \right)  \\
	& + D_z^c \left(\tens{G}_{13}  \ds  D_x^c v^n \right)
	+ D_z^c \left(\tens{G}_{23}  \ds  D_y^c v^n \right) \\
	& + D_z^- \left(A_z^+ \tens{G}_{33}  \ds  D_z^+ v^n \right).
\end{aligned}
\end{equation}
In this scheme, $D_x^-$, $D_x^+$, and $D_x^c$ denote the forward, backward, and central finite difference operators along direction $x$
while $A_x^+$ averages along direction $x$ to
provide values of a diffusion tensor component (in this case $\tens{G}_{11}$) on edge-centers on a uniform Cartesian grid.
The operator $.*$ means the component wise multiplication of arrays.
In the case of isotropic Perona--Malik diffusion, i.e., $\tens{G} = g \tens{I}$, scheme \eqref{eqn:scheme_aniso} will reduce to
the one we used in~\cite{bvgmm13}.

\subsection{Non-Dimensionalization and Parameter Adaption}
In 2D image processing one typically works on a pixel domain of the form $[0,N-1] \times [0,M-1]$ and a mesh width of $\dx = 1$.
Because our surfaces are contained in a reference box of typical size one, we work with mesh widths $\dx \ll 1$.
Thus, we are in a different scaling regime.  
In order to compare parameter choices, we non-dimensionalize the diffusion PDEs
in a standard fashion (e.g., \cite{Holmes2009foundations})
and choose parameters relative to the data.

For the linear Gaussian diffusion this results in choosing a new time-scale.
To see this, let $u$ satisfy
\begin{align*}
	\pd_t u &= \lap u, \qquad x \in [0,L_1] \times [0, L_2] \times [0,L_3],
\end{align*}
and we look at the transformed solution
\begin{align}\label{eqn:trafo}
	w(\tau,\xi) &= a + b u(\beta \tau, c + L \xi), & L &= \max\limits_{i=1,2,3}\{L_i\},
\end{align}
where we have taken into account affine linear transformations of both the color space and the domain.
The new function $w$ will then satisfy
$\pd_\tau w = \frac{\beta}{L^2} \lap_\xi w$,
and by picking the time scale $\beta = L^2$, we can then solve
$\pd_\tau w = \lap_\xi w$,
that is, the same equation for $w$ as we had for $u$, but to a different stop time.
If we are interested in $u$ at time $T$, we have to evolve $w$ until stop time $T/L^2$. 

In order to get the same scaling behavior for the non-linear Perona--Malik and Weickert models,
we adapt the parameters to the initial data (under the assumption that the initial data is non-constant).
Starting with the Perona--Malik model
\begin{align*}
	\pd_t u &= \Div \left( g(|\nabla u|^2) \nabla u \right), &
	g(s^2) = \frac{1}{1 + s^2/\lambda^2},
\end{align*}
we non-dimensionalize the nonlinear diffusivity $g$ by taking
into account the initial data $u^0$ and set
\begin{align*}
	\lambda &= \lambda_\rel \, \|\nabla u_\sigma^0\|_{\infty}.
\end{align*}
That is, a gradient is considered to be large if its magnitude is above a certain percentage $\lambda_\rel$
of the gradient magnitude of the given data $u^0$ (or precisely, the linearly smoothed version $u_\sigma^0$).
Let $w$ be again the transformation of \eqref{eqn:trafo}. Taking into account the appropriate time-scale for the linear Gaussian pre-smoothing, the data transforms as
\begin{align*}
	w^0(\xi) &= a + b u^0 (L \xi), &
	w_{\frac{\sigma}{L^2}}^0(\xi) &= a + b u_{\sigma}^0 (L \xi).
\end{align*}
Consequently, the ratio 
\begin{align*}
	\frac{|\nabla_\xi w|}{\|\nabla_\xi w_{\frac{\sigma}{L^2}}^0\|_{\infty}} = \frac{|\nabla u|}{\|\nabla u_\sigma^0\|_{\infty}}
\end{align*}
is invariant and hence we obtain
\begin{align*}
	\pd_\tau w &= \frac{\beta}{L^2} \Div_\xi \left( g(|\nabla_\xi w|^2) \nabla_\xi w \right).
\end{align*} 
By using the time-scale $\beta = L^2$ (as was the case for Gaussian diffusion) we end up with the original Perona--Malik model.

In Weickert's anisotropic diffusion models the non-linear functions depend on the coherence, so we choose the parameters relative to the coherence of the initial data.
For edge-enhancing and coherence-enhancing diffusion the eigenvalues
of the diffusion tensor $\tens{G}$ are given by \eqref{eqn:ee_ews} and
\eqref{eqn:ce_ews} respectively, in terms of the coherence~$c$.
Let $c^0$ denote the coherence of the initial data, we set
\begin{align}
	\lambda &= \lambda_\rel \, \|c^0\|_{\infty}, & B &= B_\rel \, \|c^0\|_{\infty}.
\end{align}
Because of the scaling behavior of Gaussian diffusion we have
\begin{align*}
	\tens{J}_{\sigma,\rho}[u] &= \frac{1}{b^2 L^2} \tens{J}_{\frac{\sigma}{\rule{0pt}{1.2ex}L^2},\frac{\rho}{\rule{0pt}{1.2ex}L^2}}[w].
\end{align*}
Thus, with appropriate time-scales for the Gaussian pre-smoothing and post-smoothing steps, the ratio
$c/\|c^0\|_{\infty}$ is invariant under affine linear transformations. We obtain the transformed PDE
\begin{align*}
	\pd_\tau w &= \frac{\beta}{L^2} \Div_\xi \left( \tens{G}[w] \nabla_\xi w \right).
\end{align*}
The time scale factor is thus again $\beta=L^2$.

By virtue of \eqref{eqn:sgraddiv}, the surface-intrinsic diffusion PDEs have the same scaling behavior
as their $\Real^3$-counterparts.

\section{Experiments}\label{sec:results}

We demonstrate our algorithm on several examples.
The numerical parameters used are the mesh width $\dx = 0.0125$ and the time-step size $\dt = 0.15 \dx^2$.
We use a uniform Cartesian grid defined on the reference box $[-1.5,1.5]^3$, 
but the algorithm is executed on a narrow computational band containing
only grid points close to the surface.

\subsection{Denoising}
In Figures~\ref{fig:denoising} and \ref{fig:denoising2} we compare Gaussian, Perona--Malik, and edge-enhancing diffusion.

\subsubsection{Stripe Pattern on a Torus}
\begin{figure}
	\centering
	\includegraphics[width=.48\columnwidth]{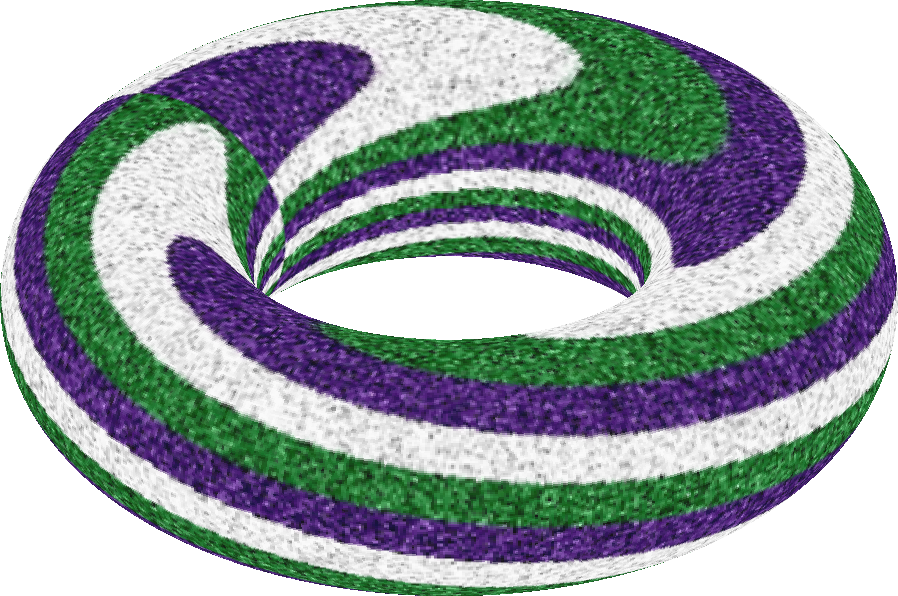}\;\;%
	\includegraphics[width=.48\columnwidth]{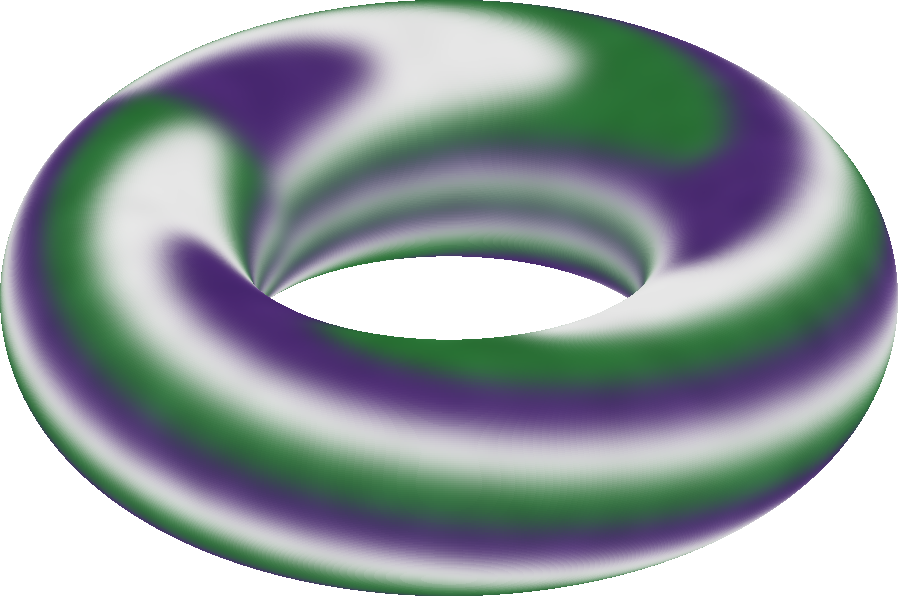}\vspace{1ex}\\
	\includegraphics[width=.48\columnwidth]{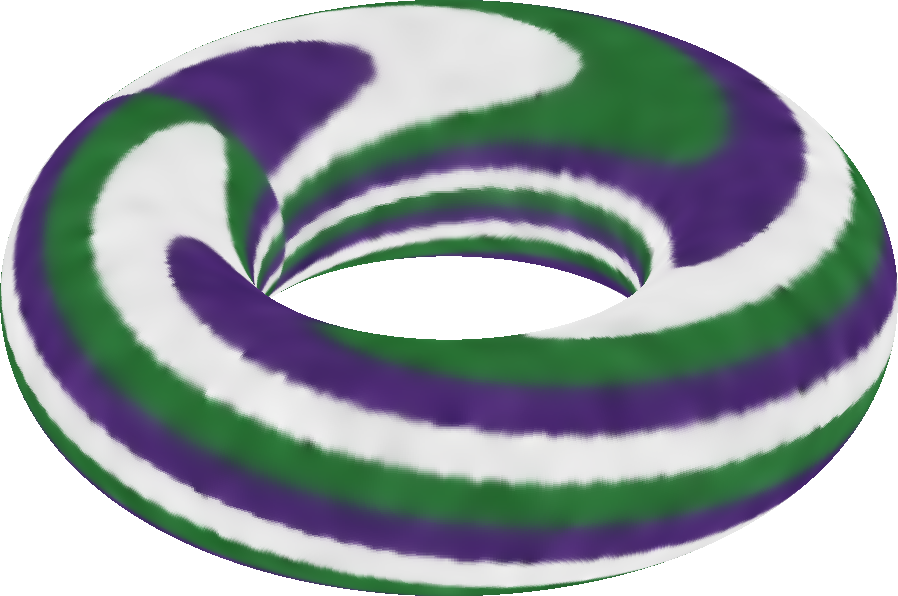}\;\;%
	\includegraphics[width=.48\columnwidth]{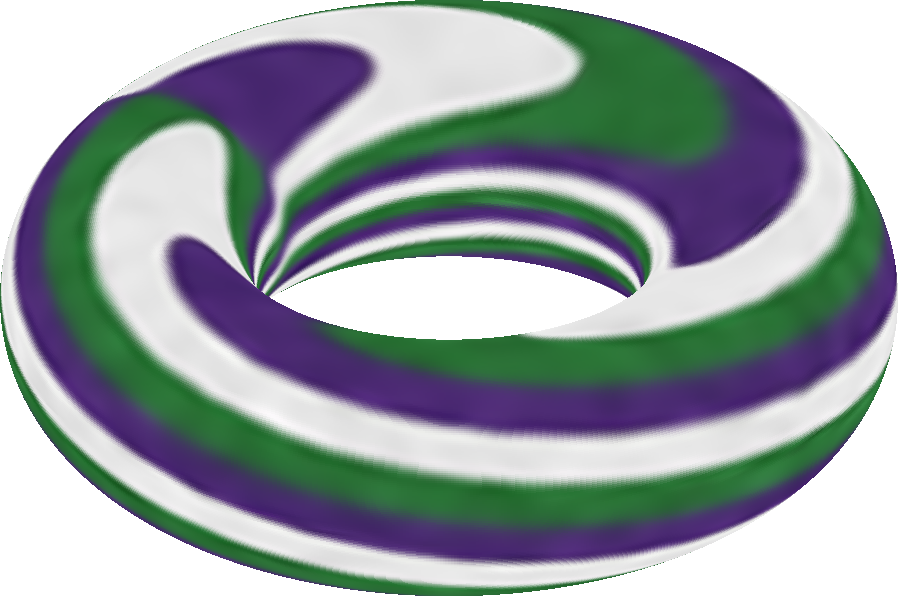}%
	\caption{Denoising on a torus: the effect of different diffusion filters.
				\textit{Top Left:} Noisy image.
				\textit{Top Right:} Gaussian diffusion.
				\textit{Bottom Left:} Perona--Malik diffusion; edge-sensitivity parameter $\lambda_\rel = 2 \cdot 10^{-1}$.
				\textit{Bottom Right:} Edge-enhancing diffusion; $\sigma = 1 \cdot 10^{-4}$, $\rho = 4 \cdot 10^{-4}$, $\lambda_\rel = 4 \cdot 10^{-2}$.
				Stop time in all cases $T=1.2\cdot 10^{-3}$ ($52$ iterations).}
	\label{fig:denoising}
\end{figure}

The top left image of Figures~\ref{fig:denoising} shows a noisy stripe pattern defined on a torus.
The torus is the solution of
\begin{align}\label{eqn:torus}
	(R - \sqrt{x^2 + y^2})^2 + z^2 &= r^2
\end{align}
where the big and the small radius are $R=1$ and $r=0.4$. 
From \eqref{eqn:torus} we find the normal field analytically and obtain the tangent space basis,
which is needed in edge-enhancing diffusion, by QR-decomposition of the normal.

In Figure~\ref{fig:denoising} we observe the usual effect of Gaussian diffusion: edges are not sharp because of the uniform blur.
Perona--Malik diffusion produces sharp edges, but the noise on the edge is still present
since the flux vanishes on edges. In contrast to that, edge-enhancing diffusion aligns the flux with edges and produces sharp but smoother edges. 

\subsubsection{Wood Grain on a Sphere} 
\begin{figure}
	\centering
	\includegraphics[width=.48\columnwidth]{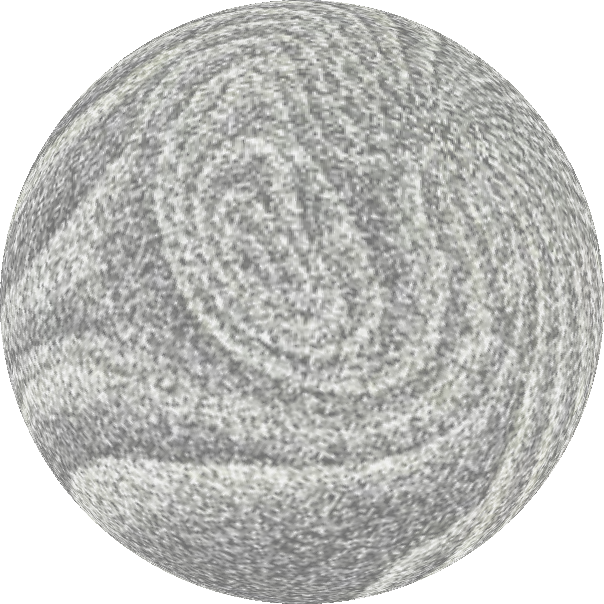}\;\;%
	\includegraphics[width=.48\columnwidth]{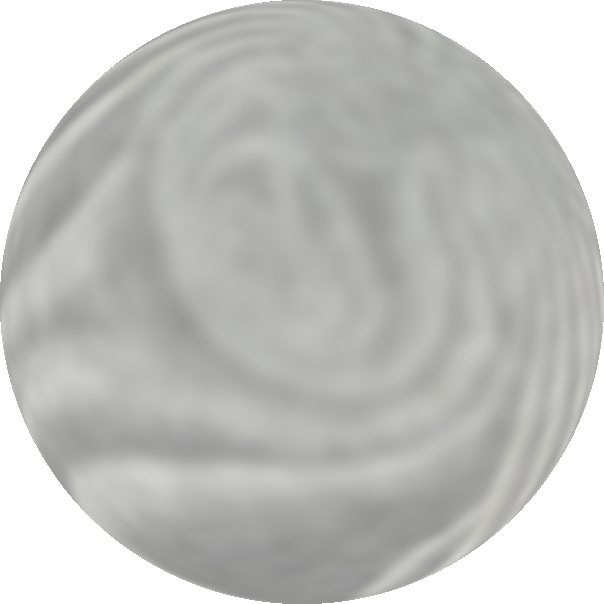}\vspace{1ex}\\
	\includegraphics[width=.48\columnwidth]{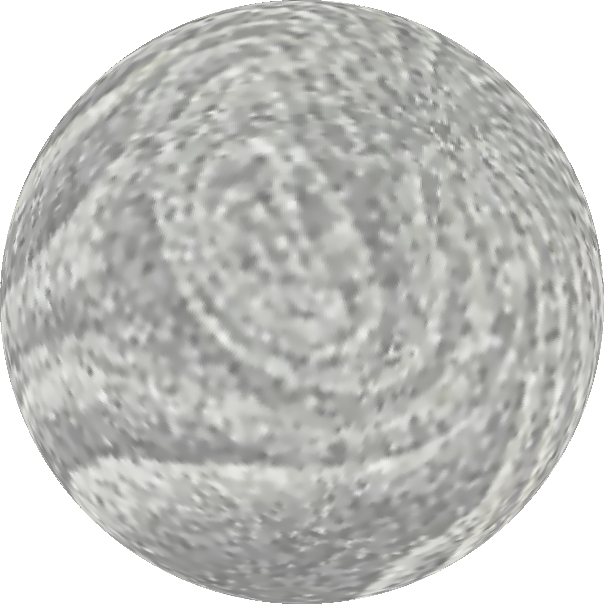}\;\;%
	\includegraphics[width=.48\columnwidth]{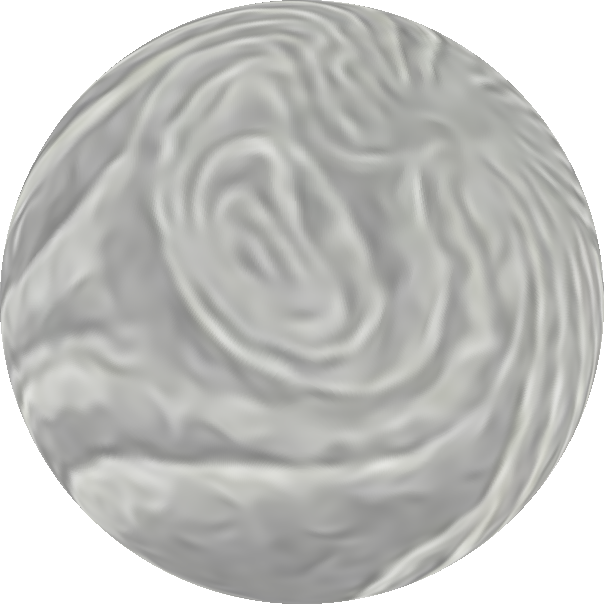}%
	\caption{Denoising on a sphere: the effect of different diffusion filters.
				\textit{Top Left:} Noisy image.
				\textit{Top Right:} Gaussian diffusion.
				\textit{Bottom Left:} Perona--Malik diffusion; edge-sensitivity parameter $\lambda_\rel = 2 \cdot 10^{-1}$.
				\textit{Bottom Right:} Edge-enhancing diffusion; $\sigma = 1 \cdot 10^{-4}$, $\rho = 4 \cdot 10^{-4}$, $\lambda_\rel = 4 \cdot 10^{-2}$.
				Stop time in all cases $T=5.9\cdot 10^{-4}$ ($25$ iterations).}
	\label{fig:denoising2}
\end{figure}

The top left image of Figures~\ref{fig:denoising2} shows a wood grain defined on the unit sphere.
The normal field is found analytically from the defining equation $x^2 + y^2 + z^2 =1$ and the tangent space basis is obtained
by QR-decomposition of the normal.

In Figure~\ref{fig:denoising2} the effects are even more drastic than in Figure~\ref{fig:denoising}:
Gaussian diffusion removes the noise but does hardly preserve any structure.
Perona--Malik diffusion preserves some structures, but most of the noise is still present.
In contrast to that, edge-enhancing diffusion removes the noise while preserving the structure of the wood grain. 

\subsection{Coherence Enhancement}

\subsubsection{Fingerprint on a Sphere}
\begin{figure}
	\centering
	\includegraphics[width=.48\columnwidth]{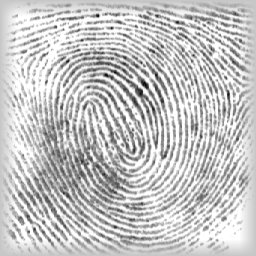}\;\;%
	\includegraphics[width=.48\columnwidth]{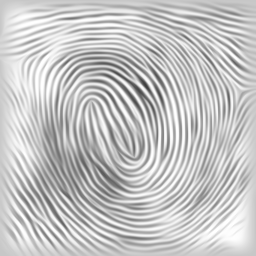}\vspace{1ex}\\
	\includegraphics[width=.48\columnwidth]{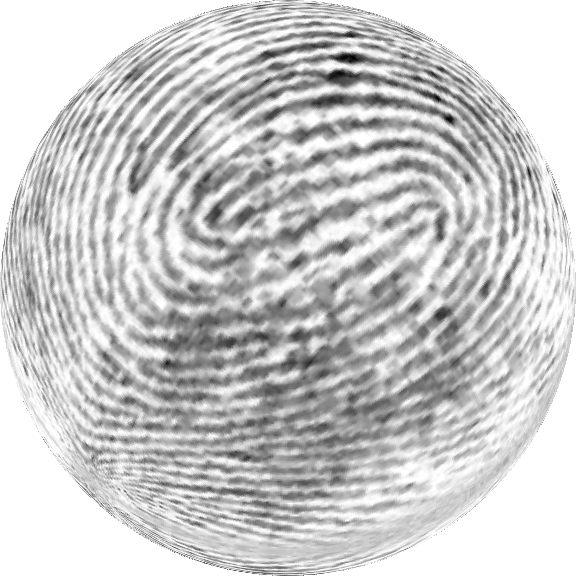}\;\;%
	\includegraphics[width=.48\columnwidth]{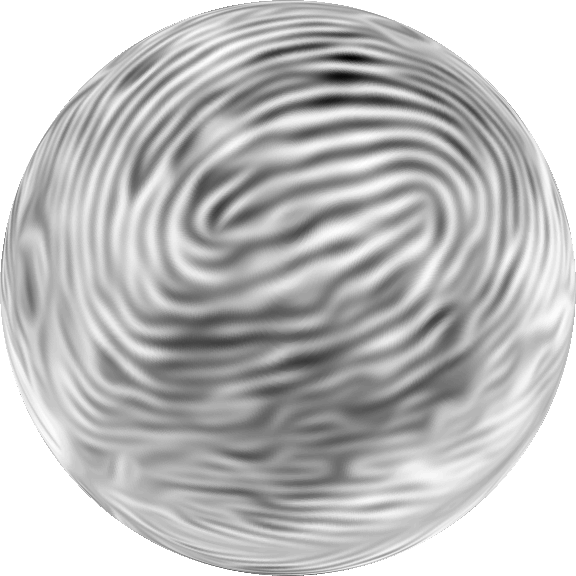}
	\caption{Feature enhancement in a fingerprint.  
				\textit{Top Left:} Fingerprint image (by courtesy of J. Weickert).
				\textit{Top Right:} Weickert's result with coherence enhancing diffusion \cite{w98}.
				\textit{Bottom Left:} Fingerprint texture-mapped to a sphere.
				\textit{Bottom Right:} Our result with in-surface coherence enhancing diffusion with the parameter values of \eqref{eqn:params_fp}.}
	\label{fig:fp}
\end{figure}

The fingerprint image is taken from Weickert \cite{w98}, here we texture-mapped it to the lower and upper hemisphere
of a unit sphere (perhaps with the upper sphere serving as simple model of finger tip).
Since the unit sphere is the solution of
$x^2 + y^2 + z^2 = 1$
we find the normal field analytically and obtain the tangent space basis by QR-decomposition of the normal.
Our result, bottom right image of Figure~\ref{fig:fp}, is similar to Weickert's result, top right image of Figure~\ref{fig:fp}.
In the following, we compare parameter choices.

Weickert reports the following parameters for his experiment with the fingerprint (an image of size $256 \times 256$ pixels with values in $[0,255]$): 
$\alpha = 10^{-3}$, $B = 1$ for the diffusion tensor and $T=20$ as the stop time for the anisotropic diffusion.
Regarding the structure tensor $\sigma = 0.5$, $\rho = 4$ were used. 

Our upper hemisphere has surface area $2\pi$ so a reasonable length factor based on area is $L = 255/\sqrt{2\pi} \approx 10^2$.
In the set-up of the structure tensor, Weickert convolves with Gaussian kernels while we are (formally) convolving with heat-kernels. 
Weickert's time scales would transform to
\begin{align*}
	\hat{\sigma} &= \frac{\sigma^2}{2 L^2} \approx 1.2 \cdot 10^{-5}, &
	\hat{\rho} &= \frac{\rho^2}{2 L^2} \approx 7.7 \cdot 10^{-4}, \\ \hat{T} &= \frac{T}{L^2} \approx 1.9 \cdot 10^{-3}.
\end{align*}
Finally we take into account our relative choice of $B$ as discussed above.
The maximal coherence of the initial state (top left image of Figure~\ref{fig:fp}) is about $\|c^0\|_\infty \approx 10^{3}$. Thus,
in order to obtain $B=1$, we have $B_\rel \approx 10^{-3}$.

The parameters which we actually used in our algorithm are
\begin{align}\label{eqn:params_fp}
	\sigma &= 10^{-4}, & \rho &= 4 \cdot 10^{-4}, & T &= 1.2 \cdot 10^{-3}, \notag \\
	\alpha &= 10^{-3}, & B_\rel &= 10^{-3},
\end{align}
and are comparable (in order of magnitude) to Weickert's.

\subsubsection{Sunflowers on a Vase}
\begin{figure}
	\centering
	\includegraphics[width=.48\columnwidth]{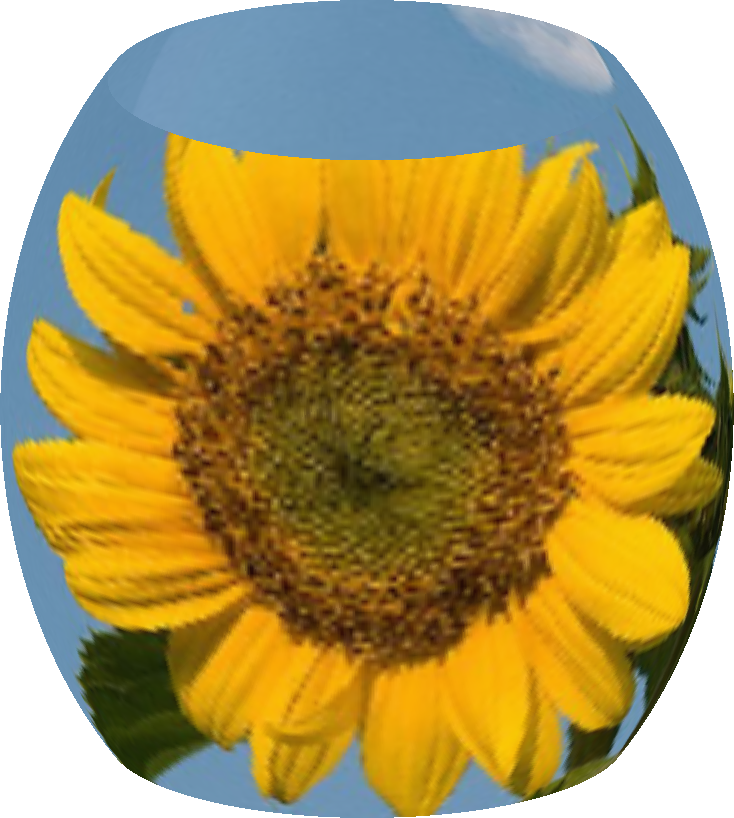}\;\;%
	\includegraphics[width=.48\columnwidth]{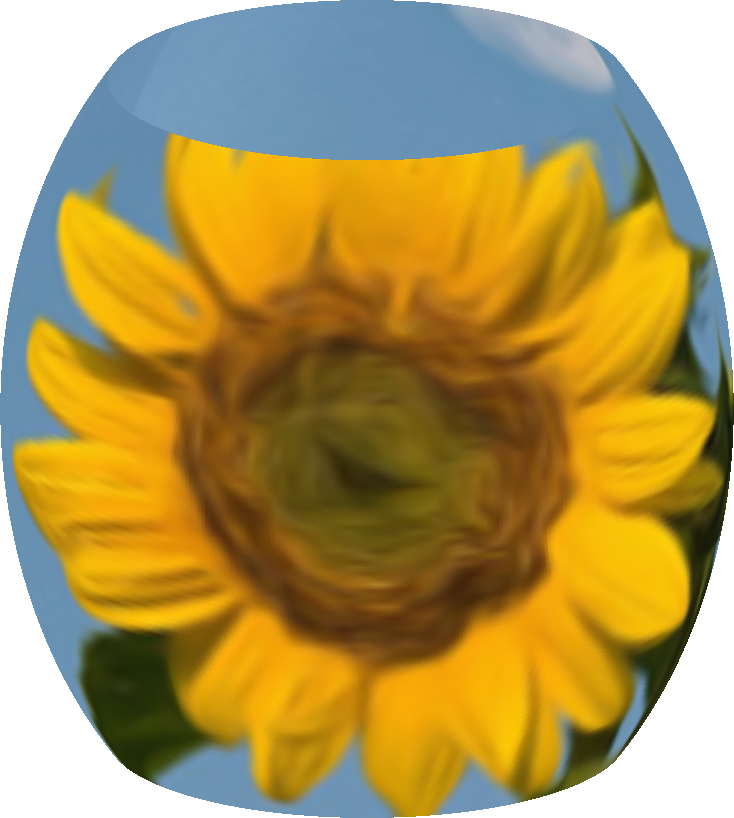}\vspace{1ex}\\
	\includegraphics[width=.48\columnwidth]{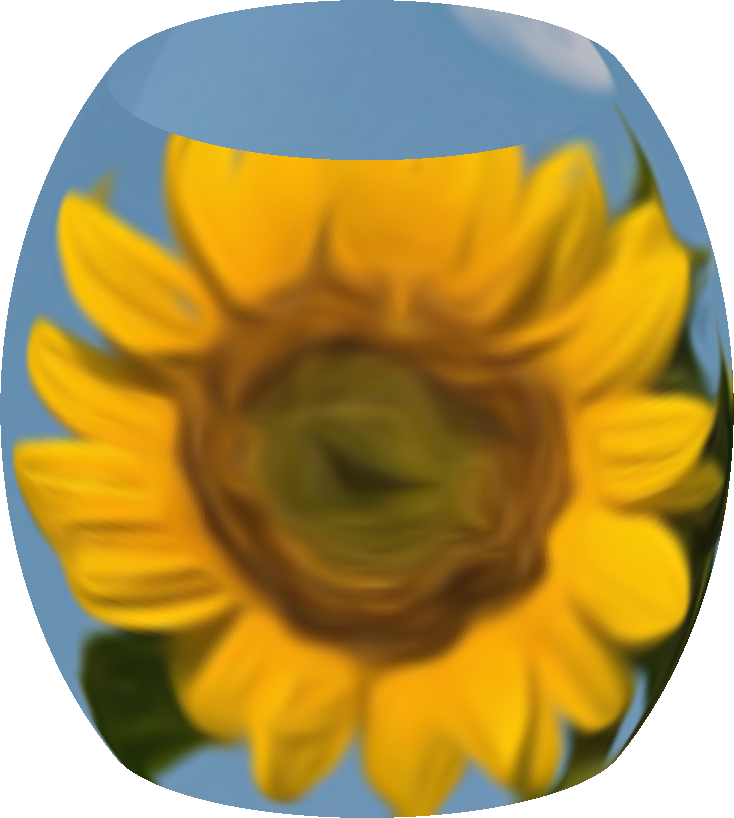}\;\;%
	\includegraphics[width=.48\columnwidth]{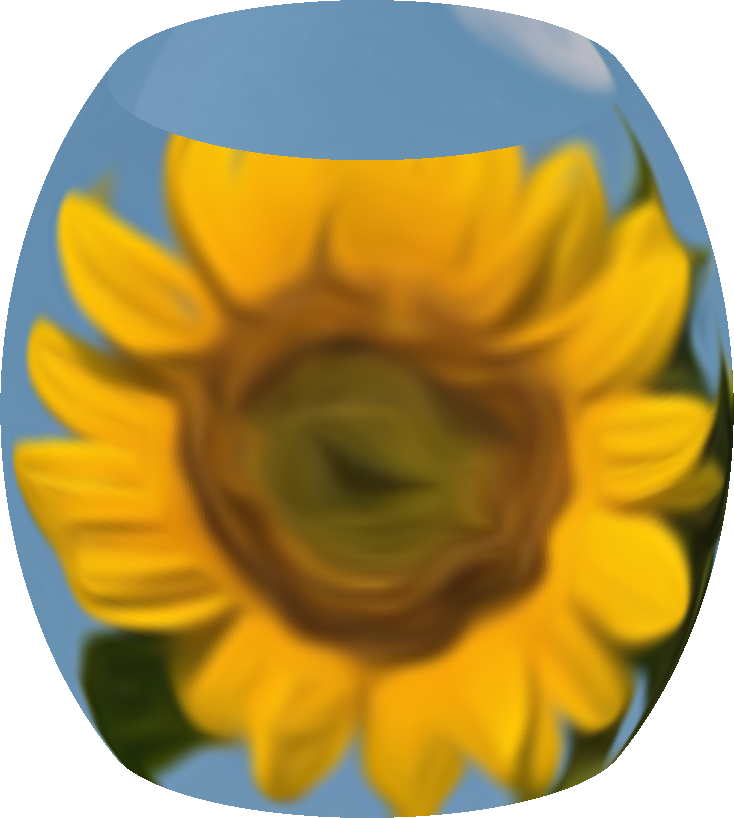}
	\caption{Stylization of a sunflower picture on a vase.
				\textit{Top Left:} Original image.
				Results of coherence enhancing diffusion with the parameter values of \eqref{eqn:params_vase}:
				\textit{Top Right:} iteration 10 ($T = 2.3 \cdot 10^{-4}$),
				\textit{Bottom Left:} iteration 20 ($T = 4.7 \cdot 10^{-4}$),
				\textit{Bottom Right:} iteration 30 ($T = 7.0 \cdot 10^{-4}$).
				}
	\label{fig:sunfl}
\end{figure}

This example demonstrates the stylization of images as an application of coherence enhancing diffusion.
The vase is a surface of revolution with a sunflower texture map.
Figure~\ref{fig:sunfl} shows the results obtained with the parameters
\begin{align}\label{eqn:params_vase}
	\sigma = 10^{-4}, \;\; \rho = 5 \cdot 10^{-4}, \;\;
	\alpha = 10^{-3}, \;\; B_\rel = 10^{-3}.
\end{align}
Since this an open surface, we must impose a boundary condition on the edge of the vase.
One reasonable choice would be the no-flux boundary condition \eqref{eqn:GSdiffBC}: $\vec{\nu}^T\vec{j} = 0$ at the boundary of $\surf$.
Another possible choice is the Neumann-zero condition $\vec{\nu}^T
\sgrad u = 0$, which we use here because it is easy to implement with
the closest point method (as it is automatically satisfied in the
closest point framework).  However, it is likely not conservative.
Nonetheless, the results of Figure~\ref{fig:sunfl} are visually reasonable at the edges.

\subsubsection{Impressionist Style on a Triangulated Surface}

\begin{figure}
  \centerline{%
    \includegraphics[width=.5\columnwidth]{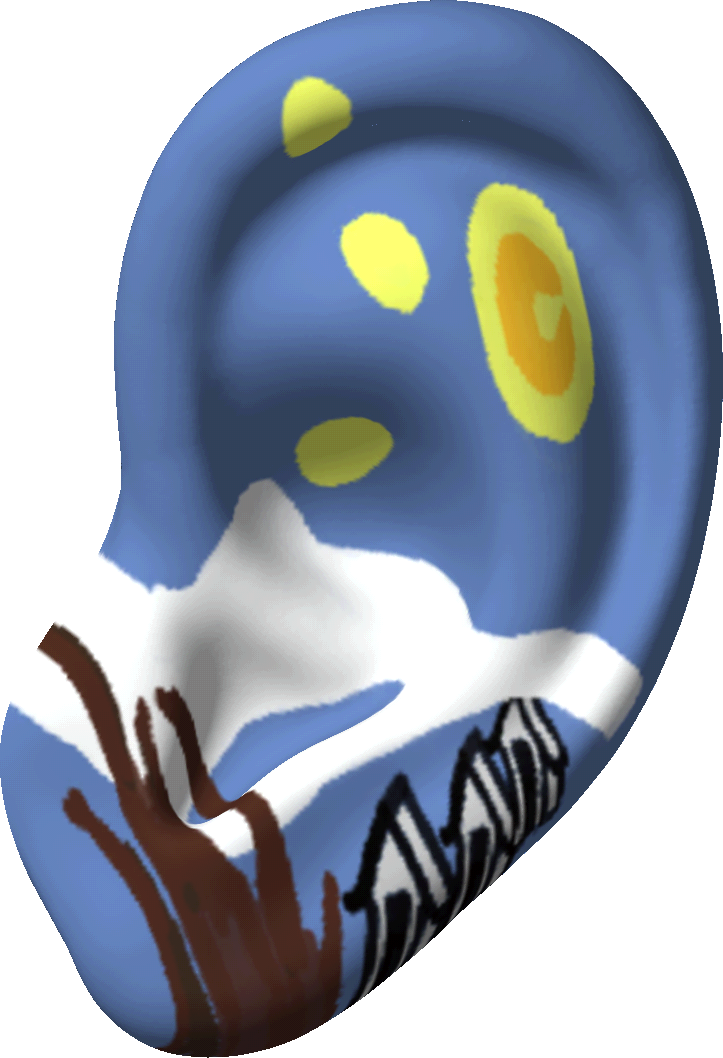}%
    \includegraphics[width=.5\columnwidth]{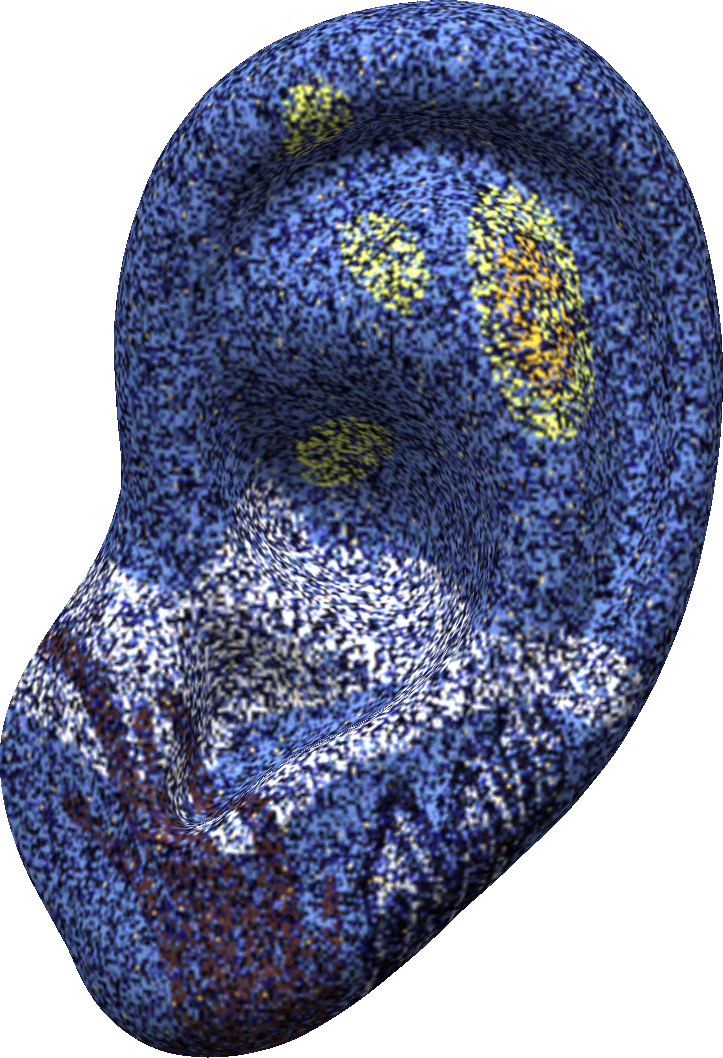}}%
  \centerline{%
    \includegraphics[width=.58\columnwidth]{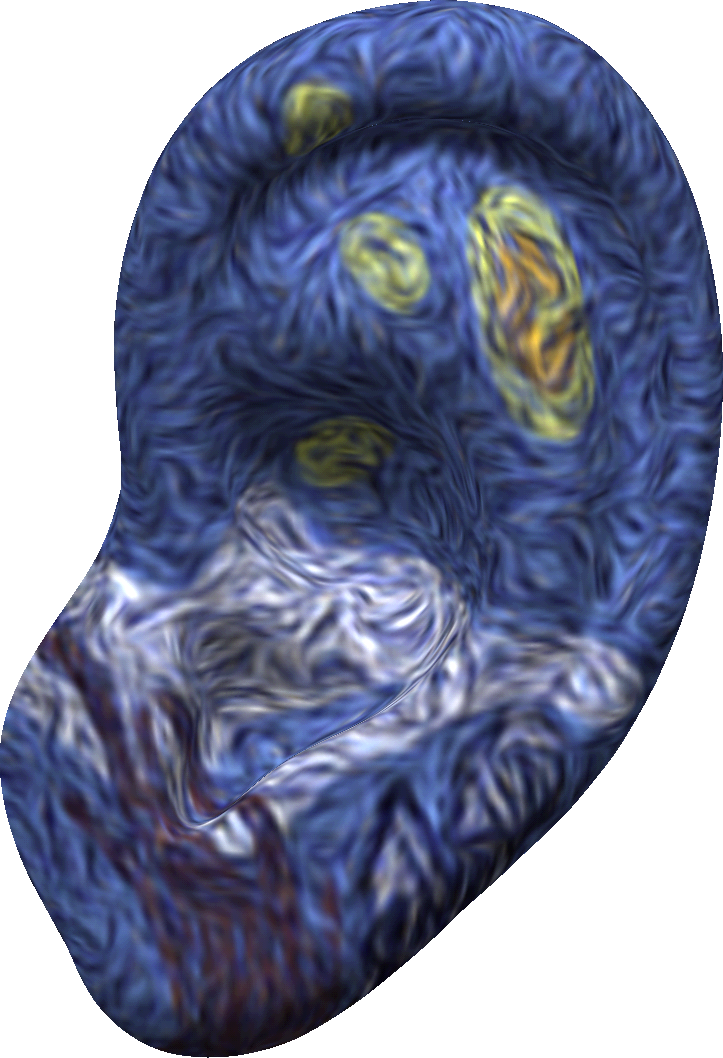}}%
  \caption{``Starry Night'' on an ear: creating Impressionist textures.
    \textit{Top-left:} user input: some colors roughly painted onto the surface,
    \textit{Top-right:} 67\% of pixels replaced with dark blue shades,
    \textit{Bottom:} Coherence enhancement, iteration 10 ($T = 1.1 \cdot 10^{-4}$).
  }
  \label{fig:ear}
\end{figure}

As another example of stylization of images using coherence enhancing
diffusion, we consider texture creation.
Here we have a triangulation of an ear, originally from
\cite{HumanEarSurface} and smoothed and Loop subdivided to a smooth mesh of
roughly 400\,000 triangles.
We convert the triangulation to a closest point function using the
technique outlined in \cite{cbm:lscpm}.
The basis for the tangent space is computed from the closest point
function using the derivative of the projection as in \eqref{eqn:Dcp}.
In Figure~\ref{fig:ear} (top-left), an artist (ahem) has quickly and
roughly painted some parts of the surface using colors chosen from the
palette of van Gogh's ``Starry Night''.
We perturb this input image by setting randomly 67\% of
the voxels to a random selection of mostly dark blues, again from the
palette of ``Starry Night''.
This noisy image is shown in Figure~\ref{fig:ear} (top-right).
The ear is two units long in its longest dimension and we embed it in
a grid with $h = 0.008$ (smaller than the other examples to get a
finer texture).
The other parameters are chosen as
\begin{align*}
	\sigma = 10^{-4}, \;\; \rho = 5 \cdot 10^{-4}, \;\;
	\alpha = 10^{-3}, \;\; B_\rel = 10^{-6}.
\end{align*}
Figure~\ref{fig:ear} (bottom) shows the results, a simulated
Impressionist painting on the surface of an ear.

\section{Conclusions} \label{sec:conclusion}
In this paper we introduced a model for edge- and coherence-enhancing image processing on curved surfaces using surface-intrinsic anisotropic diffusion.
We defined a surface-intrinsic structure tensor, from which
the construction of the diffusion tensor follows almost exactly
the procedure suggested by Weickert.
The resulting surface-intrinsic diffusion PDE is solved numerically
using the closest point method, a general method for solving PDEs
posed on surfaces.
Our approach can be used for denoising of data posed on surfaces and for visual effects such as generating surface textures or stylizing existing textures.
Our results for images on surfaces are comparable to those of Weickert and parameter choices made for 2D images
can be re-used taking into account the corresponding scale factors.

For open surfaces we used zero-Neumann rather than no-flux boundary conditions.
The implementation of general no-flux boundary conditions within the closest point method is a topic for future research.

\begin{acknowledgements}
This work was supported by award KUK-C1-013-04 made by King Abdullah University of Science and Technology (KAUST).
\end{acknowledgements}


\bibliographystyle{plain}  
\bibliography{AnisoDiffusionSurf} 

\end{document}